\documentclass{ifacconf}

\usepackage{amssymb}
\usepackage{amsmath}

\usepackage{graphicx}      
\usepackage{natbib}        
\usepackage{algorithm}
\usepackage[noend]{algpseudocode}

\newcommand{\abs}[1]{\lvert#1\rvert}

\begin{document}
\begin{frontmatter}


\title{Optimal PID and Antiwindup Control Design as a Reinforcement Learning Problem} 

\thanks[footnoteinfo]{\textcopyright 2020 the authors. This work has been accepted to IFAC World Congress for publication under a Creative Commons Licence CC-BY-NC-ND}

\author[First]{Nathan P. Lawrence} 
\author[Third]{Gregory E. Stewart} 
\author[First]{Philip D. Loewen}
\author[Fourth]{Michael G. Forbes}
\author[Fourth]{Johan U. Backstrom}
\author[Second]{R. Bhushan Gopaluni}
\address[First]{Department of Mathematics, University of British Columbia, Vancouver, BC V6T 1Z2, Canada (e-mail: lawrence@math.ubc.ca, loew@math.ubc.ca).}
\address[Second]{Department of Chemical and Biological Engineering, University of British Columbia, Vancouver, BC V6T 1Z3, Canada (e-mail: bhushan.gopaluni@ubc.ca)}
\address[Third]{Department of Electrical and Computer Engineering, University of British Columbia, Vancouver, BC V6T 1Z4, Canada (e-mail: stewartg@ece.ubc.ca)}
\address[Fourth]{Honeywell Process Solutions, North Vancouver, BC V7J 3S4, Canada (e-mail: michael.forbes@honeywell.com, johan.backstrom@honeywell.com)}

\begin{abstract}                
Deep reinforcement learning (DRL) has seen several successful applications to process control. Common methods rely on a deep neural network structure to model the controller or process. With increasingly complicated control structures, the closed-loop stability of such methods becomes less clear. In this work, we focus on the interpretability of DRL control methods. In particular, we view linear fixed-structure controllers as shallow neural networks embedded in the actor-critic framework. PID controllers guide our development due to their simplicity and acceptance in industrial practice. We then consider input saturation, leading to a simple nonlinear control structure. In order to effectively operate within the actuator limits we then incorporate a tuning parameter for anti-windup compensation. Finally, the simplicity of the controller allows for straightforward initialization. This makes our method inherently stabilizing, both during and after training, and amenable to known operational PID gains.
\end{abstract}


\end{frontmatter}

\section{Introduction}

The performance of model-based control methods such as model predictive control (MPC) or internal model control (IMC) relies on the accuracy of the available plant model. Inevitable changes in the plant over time result in increased plant-model uncertainty and decreased performance of the controllers. Model reidentification is costly and time-consuming, often making this procedure impractical and less frequent in industrial practice.


Reinforcement learning (RL) is a branch of machine learning in which the objective is to learn an optimal policy (controller) through interactions with a stochastic environment \citep{sutton2018reinforcement}. Only somewhat recently has RL been successfully applied in the process industry \citep{badgwell2018reinforcement}. The first successful implementations of RL methods in process control were developed in the early 2000s. For example, \citet{lee2001neuro, lee2008value} utilize approximate dynamic programming (ADP) methods for optimal control of discrete-time nonlinear systems. While these results illustrate the applicability of RL in controlling discrete-time nonlinear processes, they are also limited to processes for which at least a partial model is available or can be derived through system identification.



Other approaches to RL-based control use a fixed control structure such as PID.
With applications to process control, \citet{brujeni2010dynamic} develop a model-free algorithm to dynamically assign the PID gains from a pre-defined collection derived from IMC. On the other hand, \cite{berger2013neurodynamic} dynamically tune a PID controller in continuous parameter space using the actor-critic method, where the actor is the PID controller; the approach is based on dual heuristic dynamic programming, where an identified model is assumed to be available. The actor-critic method is also employed in \cite{sedighizadeh2008adaptive}, where the PID gains are the actions at each time-step.


\begin{figure}
\begin{center}
\includegraphics[width=8.4cm]{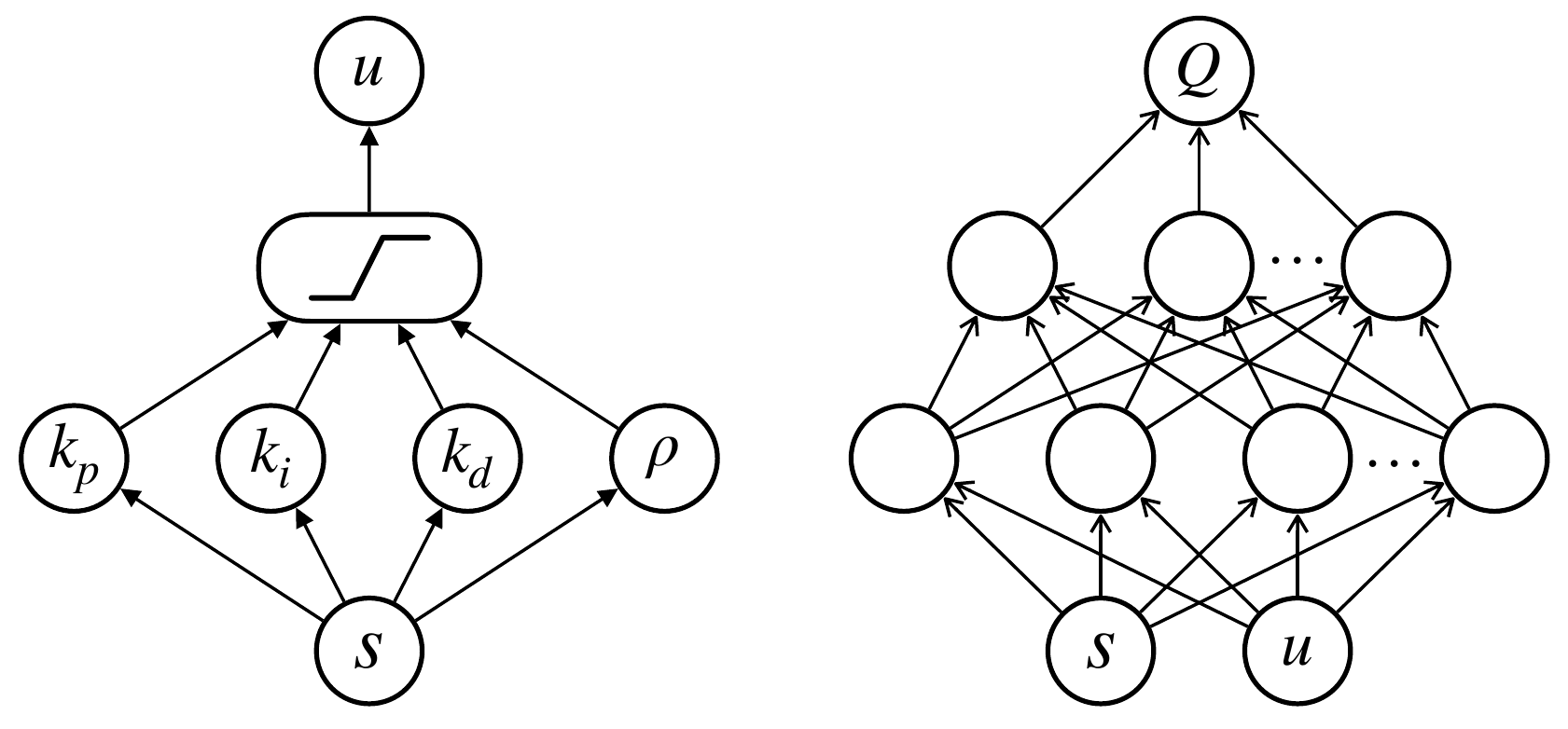}    
\caption{The actor (PID controller) on the left is a linear combination of the state and the PID \& anti-windup parameters followed by a nonlinear saturation function. The critic on the right is a deep neural network approximation of the $Q$-function whose inputs are the state-action pair generated by the actor.} 
\label{fig:PID_NN}
\end{center}
\end{figure}

In the recent work of \citet{spielberg2019toward}, an actor-critic architecture based on the deep deterministic policy gradient (DDPG) algorithm due to \cite{lillicrap2015continuous} is implemented to develop a model-free, input-output controller for set-point tracking problems of discrete-time nonlinear processes. The actor and critic are both parameterized by ReLU deep neural networks (DNNs). At the end of training, the closed-loop system includes a plant together with a neural network as the nonlinear feedback controller. The neural network controller is a black-box in terms of its stabilizing properties. In contrast, PID controllers are widely used in industry due to their simplicity and interpretability. However, PID tuning is also known to be a challenging nonlinear design problem, making it an important and practical baseline for RL algorithms.

To this end, we present a simple interpretation of the actor-critic framework by expressing a PID controller as a shallow neural network (figure \ref{fig:PID_NN} illustrates the proposed framework). The PID gains are the weights of the actor network. The critic is the $Q$-function associated with the actor, and is parameterized by a DNN. We then extend our interpretation to include input saturation, making the actor a simple nonlinear controller. Input saturation can lead to integral windup; we therefore incorporate a new tuning parameter for anti-windup compensation.  Finally, the simplicity of the actor network allows us to initialize training with hand-picked PID gains, for example, with SIMC \citep{skogestad2001probably}. The actor is therefore initialized as an operational, interpretable, and industrially accepted controller that is then updated in an optimal direction after each roll-out (episode) in the plant. Although a PID controller is used here, the interpretation as a shallow neural network applies for any linear fixed-structure controller.


This paper is organized as follows: Section \ref{sec:PID} provides a brief description of PID control and anti-windup compensation. Section \ref{sec:RL} frames PID tuning in the actor-critic architecture and describes our methodology and algorithm. Finally, section \ref{sec:results} shows simulation results in tuning a PI controller as well as a PI controller with anti-windup compensation.

\section{PID control and Integral Windup}\label{sec:PID}

We use the parallel form of the PID controller:
\begin{equation}\label{eq:PIDt}
    u(t) = k_p e_y(t) + k_i \int_{0}^{t} e_y(\tau) d\tau + k_d \frac{d}{dt}e_y(t).
\end{equation}
Here $e_{y}(t) = \bar{y}(t) - y(t)$ is the difference between the output $y(t)$ and a given reference signal $\bar{y}(t)$.
To implement the PID controller it is necessary to discretize in time. Let $\Delta t > 0$ be a fixed sampling time. Then define $I_{y}(t_n) = \sum_{i = 1}^{n} e_{y}(t_i) \Delta t$, where $0 = t_0 < t_1 < \ldots < t_n$, and $D(t_n) = \frac{e_y(t_n) - e_y(t_{n-1})}{\Delta t}$. We then use $u$ to refer to the discretized version of (\ref{eq:PIDt}), written as follows
\begin{equation}
u(t_n) = k_p e_y(t_n) + k_i I_{y}(t_n) + k_d D(t_n).
\label{eq:PIDd}
\end{equation}
We note that the velocity form of a PID controller could also be used in the following sections. However, it is simpler and more common to explain our anti-windup strategy with the form of equation \eqref{eq:PIDd}. Further, the velocity form is more sensitive to noise (exploration noise) added to the input because it gets carried over to subsequent time-steps.

Despite their simplicity, PID controllers are difficult to tune for a desired performance. Popular strategies for PID tuning include relay tuning (e.g., \cite{aastrom1984automatic}) and IMC (e.g., \cite{skogestad2001probably}). This difficulty can be exacerbated when a PID controller is implemented on a physical plant due to the limitations of an actuator. In the next section, we describe how such limitations can be problematic, then introduce a practical and simple method for working within these constraints.


\subsection{Anti-Windup Compensation}

A controller can become saturated when it has maximum and minimum constraints on its control signal and is given a set-point or a disturbance that carries the control signal outside these limits. If the actuator constraints are given by two scalars $u_{\text{min}} < u_{\text{max}}$, then we define the saturation function to be
\begin{equation}
\text{sat}(u) = 
\begin{cases}
u_{\text{min}}, &\text{ if } u < u_{\text{min}}\\
u,  &\text{ if } u_{\text{min}} \leq u \leq u_{\text{max}}\\
u_{\text{max}}, &\text{ if } u > u_{\text{max}}.
\end{cases}
\label{eq:sat}
\end{equation}
If saturation persists, the controller is then operating in open-loop and the integrator continues to accumulate error at a non-diminishing rate. That is, the integrator experiences \emph{windup}. This creates a nonlinearity in the controller and can destabilize the closed-loop system. Methods for mitigating the effects of windup are referred to as \emph{anti-windup} techniques. For a more detailed overview of the windup phenomenon and simple anti-windup techniques, the reader is referred to \cite{astrom1989integrator}.

\begin{figure}[tb]
\begin{center}
\includegraphics[width=8.4cm]{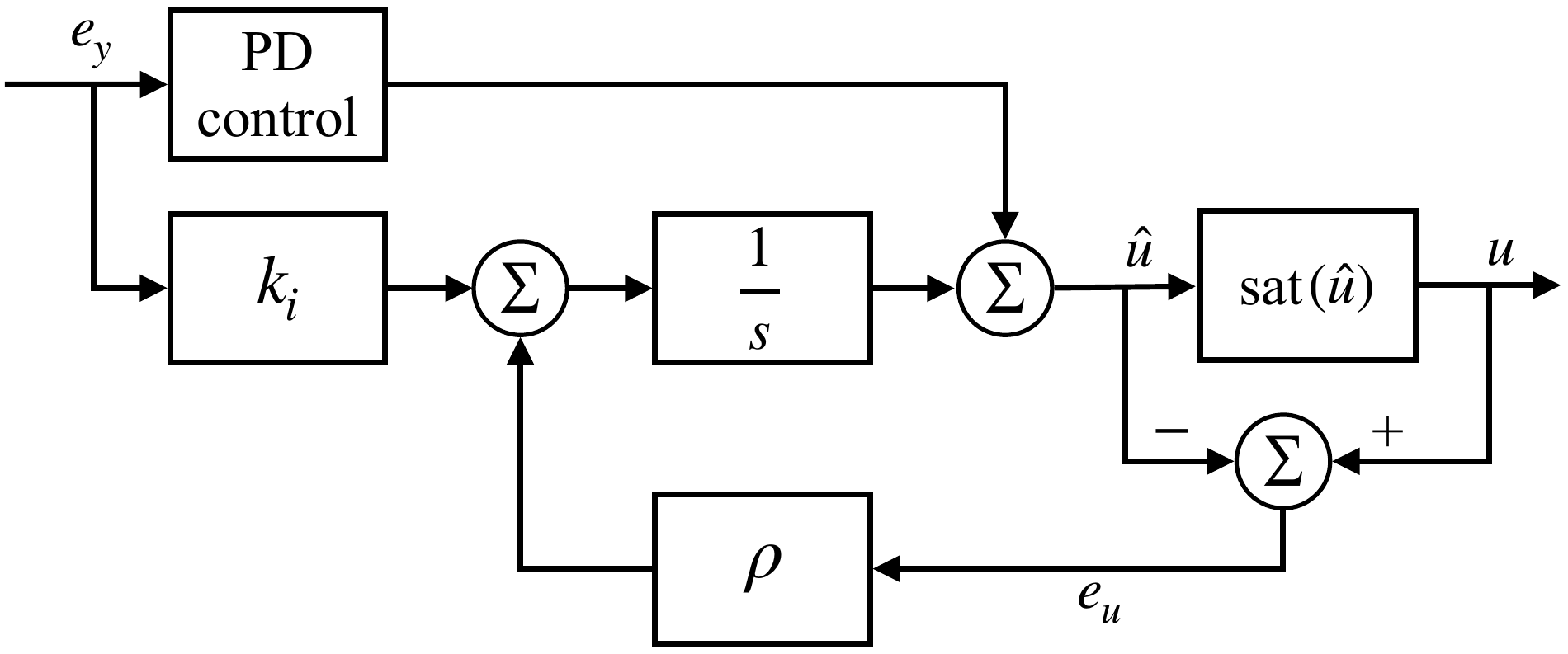}    
\caption{The back-calculation scheme feeds the scaled difference between the saturated input signal and that suggested by a PID controller back into the integrator.} 
\label{fig:Antiwindup}
\end{center}
\end{figure}

In this paper, we focus on one of the earliest and most basic anti-windup methods called \emph{back-calculation} \citep{fertik1967direct}. Back-calculation works in discrete-time by feeding into the control signal a scaled sum of past deviations of the actuator signal from the unsaturated signal. The nonnegative scaling constant, $\rho$, governs how quickly the controller unsaturates (that is, returns to the interval $(u_{\text{min}}, u_{\text{max}})$). Precisely, we define $e_{u}(t) = \text{sat}(u(t)) - u(t)$ and $I_{u}(t_n) = \sum_{i = 1}^{n-1} e_{u}(t_i) \Delta t$, then we redefine the PID controller in \eqref{eq:PIDd} to be the following
\begin{equation}
    u(t_n) = k_p e_y(t_n) + k_i I_{y}(t_n) + k_d D(t_n) + \rho I_{u}(t_n)
\label{eq:PIDd_AW}
\end{equation}
From \eqref{eq:sat} it is clear that if the controller is operating within its constraints, then \eqref{eq:PIDd_AW} emits the same control signal as \eqref{eq:PIDd}. Otherwise, the difference $\text{sat}(u) - u$ adds negative feedback to the controller if $u > u_{\text{max}}$, or positive feedback if $u < u_{\text{min}}$. Further, \eqref{eq:PIDd_AW} agrees with \eqref{eq:PIDd} when $\rho = 0$; therefore, the recovery time of the controller to the operating region $[u_{\text{min}}, u_{\text{max}}]$ is slower the closer $\rho$ is to zero and more aggressive when $\rho$ is large. A scheme of this approach is shown in figure \ref{fig:Antiwindup} and the effect of the parameter $\rho$ is shown in figure \ref{fig:ex3_fig1} in Section \ref{sec:results}.

\section{PID in the Reinforcement Learning Framework}\label{sec:RL}

Our method for PID tuning stems from the state-space representation of (\ref{eq:PIDd_AW}) followed by the input saturation:

\begin{align}
\begin{split}
\label{eq:SS1}
    \begin{bmatrix}
    e_{y}(t_{n})\\
    I_{y}(t_{n})\\
    D(t_{n})\\
    I_{u}(t_{n})
    \end{bmatrix}
    =
    &
    \begin{bmatrix}
    0 & 0 & 0 & 0\\
    0 & 1 & 0 & 0\\
    -1/\Delta t & 0 & 0 & 0\\
    0 & 0 & 0 & 1
    \end{bmatrix}
    \begin{bmatrix}
    e_{y}(t_{n-1}) \\
    I_{y}(t_{n-1}) \\
    D(t_{n-1}) \\
    I_{u}(t_{n-1})
    \end{bmatrix}
    + \\
    &\begin{bmatrix}
    1 & 0 \\
    \Delta t & 0 \\
    1/\Delta t & 0 \\
    0 & \Delta t
    \end{bmatrix}
    \begin{bmatrix}
    e_{y}(t_{n}) \\
    e_{u}( t_{n}
    \end{bmatrix}
\end{split}\\
\label{eq:SS2}
    \hat{u}(t_{n}) =
    &
    \begin{bmatrix}
    k_p & k_i & k_d & \rho
    \end{bmatrix}
    \begin{bmatrix}
    e_{y}(t_{n}) \\
    I_{y}(t_{n}) \\
    D(t_{n}) \\
    I_{u}(t_{n})
    \end{bmatrix}\\
    u(t_{n}) =& \text{ sat}\big(\hat{u}(t_{n})\big).\label{eq:SS3}
\end{align}

Equation (\ref{eq:SS1}) describes the computations necessary for implementing a PID controller in discrete time steps. On the other hand, \eqref{eq:SS2} parameterizes the PID controller. We therefore take \eqref{eq:SS2} and \eqref{eq:SS3} to be a shallow neural network, where $[k_p\ k_i\ k_d\ \rho]$ is a vector of trainable weights and the saturation is a nonlinear activation.

In the next section we outline how RL can be used to train these weights without a process model. The overview of RL provided here is brief. For a thorough tutorial of RL in the context of process control the reader is referred to the paper by \citet{spielberg2019toward}. Further, the general DDPG algorithm we employ is introduced by \citet{lillicrap2015continuous}.

\subsection{Overview of Tuning Objective}\label{subsec:objective}

The fundamental components of RL are the policy, the objective, and the environment. We assume the environment is modeled by a Markov decision process with action space $\mathcal{U}$ and state space $\mathcal{S}$. Therefore, the environment is modeled with an initial distribution $p(s_0)$ with a transition distribution $p(s_{n+1} | s_n, u_n )$, where $s_0, s_{n}, s_{n+1} \in \mathcal{S}$ and $u_{n} \in \mathcal{U}$. Here, we define $s_n = [e_{y}(t_n)\ I_{y}(t_n)\ D(t_n)\ I_{u}(t_n)]^T$ and $u_n \in \mathcal{U}$ refers to the saturated input signal given by \eqref{eq:SS3} at time $t_n$. The vector of parameters in \eqref{eq:SS2} is referred to as $K$. Formally, the PID controller with anti-windup compensation in \eqref{eq:SS3} is given by the mapping
$\mu( \cdot, K )\colon\mathcal{S}\rightarrow\mathcal{U}$ such that
\begin{equation}
u_n = \mu(s_n, K)
\label{eq:actor}
\end{equation}

Each interaction the controller \eqref{eq:actor} has with the environment is scored with a scalar value called the \emph{reward}. Reward is given by a function $r: \mathcal{S} \times \mathcal{U} \to \mathbb{R}$; we use $r_n$ to refer to $r(s_n, u_n)$ when the corresponding state-action pair is clear. We use the notation $h\sim p^{\mu}(\cdot)$ to denote a trajectory $h = (s_1, u_1, r_1, \ldots, s_N, u_N, r_N)$ generated by the policy $\mu$, where $N$ is a random variable called the \emph{terminal time}.

The goal of RL is to find a controller, namely, the weights $K$, that maximizes the expectation of future rewards over trajectories $h$:
\begin{equation}
J(\mu( \cdot, K) ) 
= \mathbb{E}_{h\sim p^{\mu}(\cdot)}\bigg[\sum_{n=1}^{\infty} \gamma^{n-1}r(s_n,\mu(s_n, K))\bigg| s_0 \bigg]
\label{eq:DPG_Obj}
\end{equation}
where $s_0 \in \mathcal{S}$ is a starting state and $0 \leq \gamma \leq 1$ is a \emph{discount} factor. Our strategy is to iteratively maximize $J$ via stochastic gradient ascent, as maximizing $J$ corresponds to finding the optimal PID gains. Optimizing this objective requires additional concepts, which we outline in the next section.

\subsection{Controller Improvement}

Equation \eqref{eq:DPG_Obj} is referred to as the \emph{value function} for policy $\mu$.  Closely related to the cost function is the \emph{Q-function}, or \emph{state-action value function}, which considers state-action pairs in the conditional expectation:
\begin{equation}
    Q(s_n, u_n) = \mathbb{E}_{h\sim p^{\mu}(\cdot)}\bigg[\sum_{k=n}^{\infty} \gamma^{k-n}r(s_k,\mu(s_k, K))\bigg| s_n, u_n \bigg]\\
\label{eq:Q_fun}    
\end{equation}
Returning to our objective of maximizing \eqref{eq:DPG_Obj}, we employ the policy gradient theorem for deterministic policies \citep{silver2014deterministic}:
\begin{equation}
\label{eq:DPGOffPolicy}
\begin{split}
&\nabla_{K} J(\mu(\cdot, K)) = \\
&\mathbb{E}_{h \sim 
p^{\mu}(\cdot)}\big[\nabla_u 
Q(s_n, u)|_{u=\mu(s_n, K)}\nabla_{K} \mu(s_n, K)\big].
\end{split}
\end{equation}
We note that Eq. (\ref{eq:DPG_Obj}) is maximized only when the policy parameters $K$ are optimal, 
which then leads to the update scheme
\begin{equation}
K
\leftarrow  K+\alpha\nabla_{K} J(\mu(\cdot, K)),
\label{eq:PolicyGradient_Iteration}
\end{equation}
where $\alpha>0$ is the learning rate.




\subsection{Deep Reinforcement Learning}

The optimization of $J$ in line \eqref{eq:PolicyGradient_Iteration} relies on knowledge of the $Q$-function \eqref{eq:Q_fun}. We approximate $Q$ iteratively using a deep neural network with training data from replay memory (RM). RM is a fixed-size collection of tuples of the form $(s_n, u_n, s_{n+1}, r_{n})$. Concretely, we write a parametrized $Q$-function, $Q( \cdot, \cdot, W_c): \mathcal{S} \times \mathcal{U} \to \mathbb{R}$, where $W_c$ is a collection of weights. This framework gives rise to a class of RL methods known as \emph{actor-critic} methods. Precisely, the actor-critic methods utilize ideas from policy gradient methods and $Q$-learning with function approximation \citep{konda2000actor, sutton2000policy}. Here, the actor is the PID controller given by \eqref{eq:actor} and the critic is $Q( \cdot, \cdot, W_c)$. 



\subsection{Actor-Critic Initialization}\label{subsec:initialization}

An advantage of our approach is that the weights for the actor can be initialized with user-specified PID gains. For example, if a plant is operating with known gains $k_p, k_i$, and $k_d$, then these can be used to initialize the actor. The idea is that these gains will be updated by stochastic gradient ascent in the approximate direction leading to the greatest expected reward. The quality of the gain updates then relies on the quality of the $Q$-function used in (\ref{eq:PolicyGradient_Iteration}). The $Q$-function is parameterized by a deep neural network and is therefore initialized randomly. Both the actor and critic parameters are updated after each roll-out with the environment. However, depending on the number of time-steps in each roll-out, this can lead to slow learning. Therefore, we continually update the critic during the roll-out using batch data from RM.  

\subsection{Connections to Gain Scheduling}\label{subsec:gain_schedule}

As previously described, the actor (PID controller) is updated after each episode. We are, however, free to change the PID gains at each time-step. In fact, previous approaches to RL-based PID tuning such as \cite{brujeni2010dynamic} and \cite{sedighizadeh2008adaptive} dynamically change the PID gains at each time-step. 
There are two main reasons for avoiding this whenever possible. 
First, the PID controller is designed for set-point tracking and is an inherently intelligent controller that simply needs to be improved subject to the user-defined objective (reward function); that is, it does not need to `learn' how to track a set-point.   
Second, when the PID gains are free to change at each time-step, the policy essentially functions as a gain scheduler. This switching of the control law creates nonlinearity in the closed-loop, making the stability of the overall system more difficult to analyze. This is true even if all the gains or controllers involved are stabilizing \citep{stewart2012pragmatic}. See, for instance, example 1 of \cite{malmborg1996stabilizing}.

Of course, gain scheduling is an important strategy for industrial control. The main point here is that RL-based controllers can inherit the same stability complications as gain scheduling. In the next subsection, we demonstrate the effect of updating the actor at different rates on a simple linear system.

\section{Simulation Results}\label{sec:results}

In our examples we refer to several different versions of the DDPG algorithm which are differentiated based on how frequently the actor is updated: \textbf{V1} updates the actor at each time-step, while \textbf{V2} updates the actor at the end of each episode. See Appendix \ref{app:details} for implementation details.

For our purposes, we define the reward function to be
\begin{equation}
    r(s_n, u_n) = -\big(\abs{e_{y}(t_n)}^p + \lambda \abs{u_n}\big),
    \label{eq:reward}
\end{equation}
where $p \in \left\lbrace 1,2\right\rbrace$ and $\lambda \geq 0$ are fixed during training. An episode ends either after 200 time-steps or when the actor tracks the set-point for 10 time-steps consecutive time-steps. 

\begin{figure}
\begin{center}
\includegraphics[width=8.4cm]{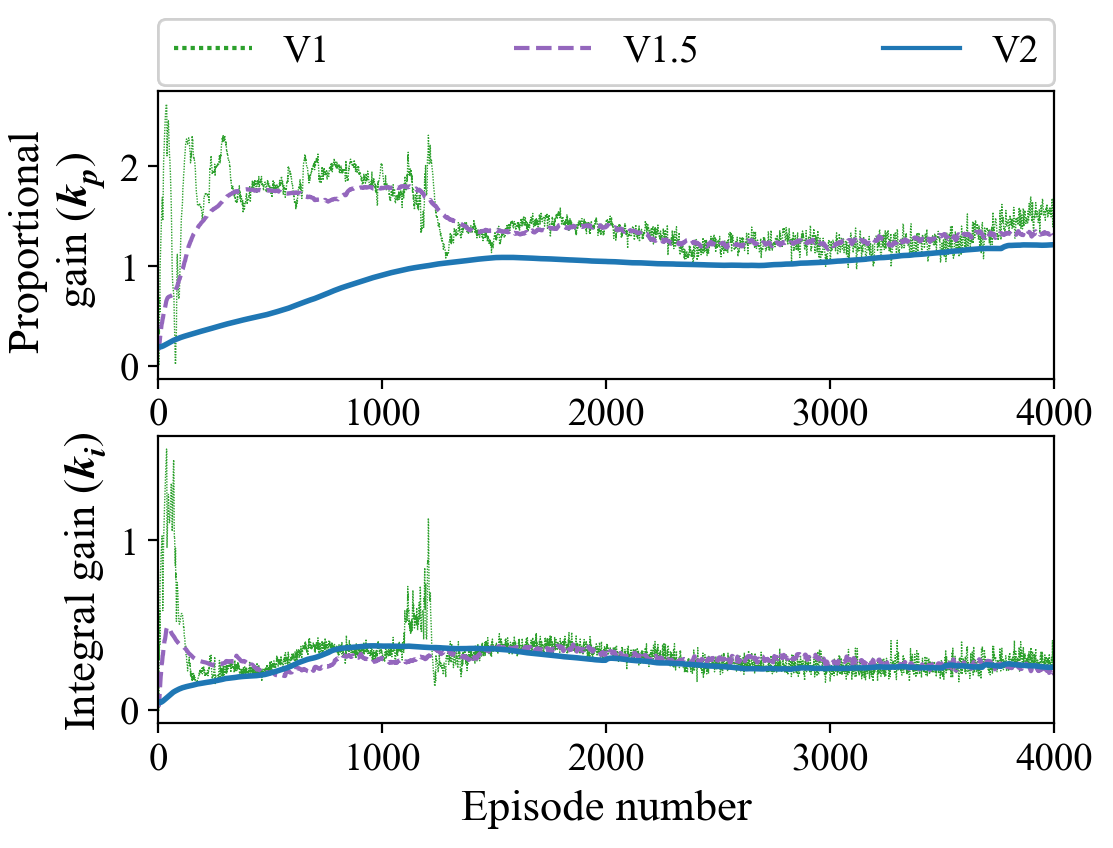}    
\caption{(top) $k_p$ parameter values at the end of each episode; (bottom) similarly, the $k_i$ parameter values. Green corresponds to updating the PI parameters at each time-step; purple corresponds to an update every 10th time-step; blue corresponds to a single update per episode. The color scheme is consistent throughout the example.} 
\label{fig:ex1_fig2}
\end{center}
\end{figure}

\begin{figure}
\begin{center}
\includegraphics[width=8.4cm]{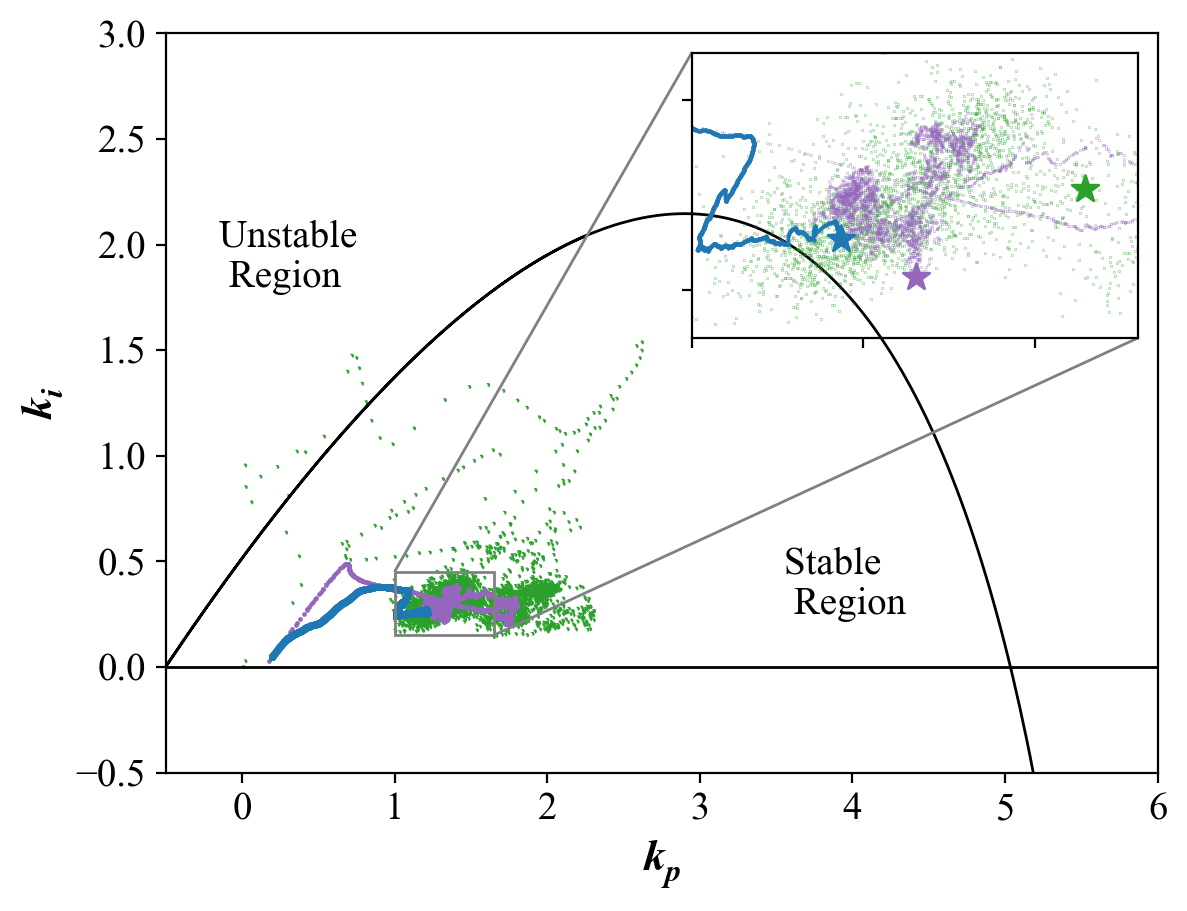}    
\caption{A scatter plot of the data shown in figure \ref{fig:ex1_fig2}. The black curve indicates the boundary of stability in the parameter plane. Stars show the $k_p-k_i$ coordinate at the end of training for its respective color.} 
\label{fig:ex1_fig5}
\end{center}
\end{figure}

\begin{figure}
\begin{center}
\includegraphics[width=8.4cm]{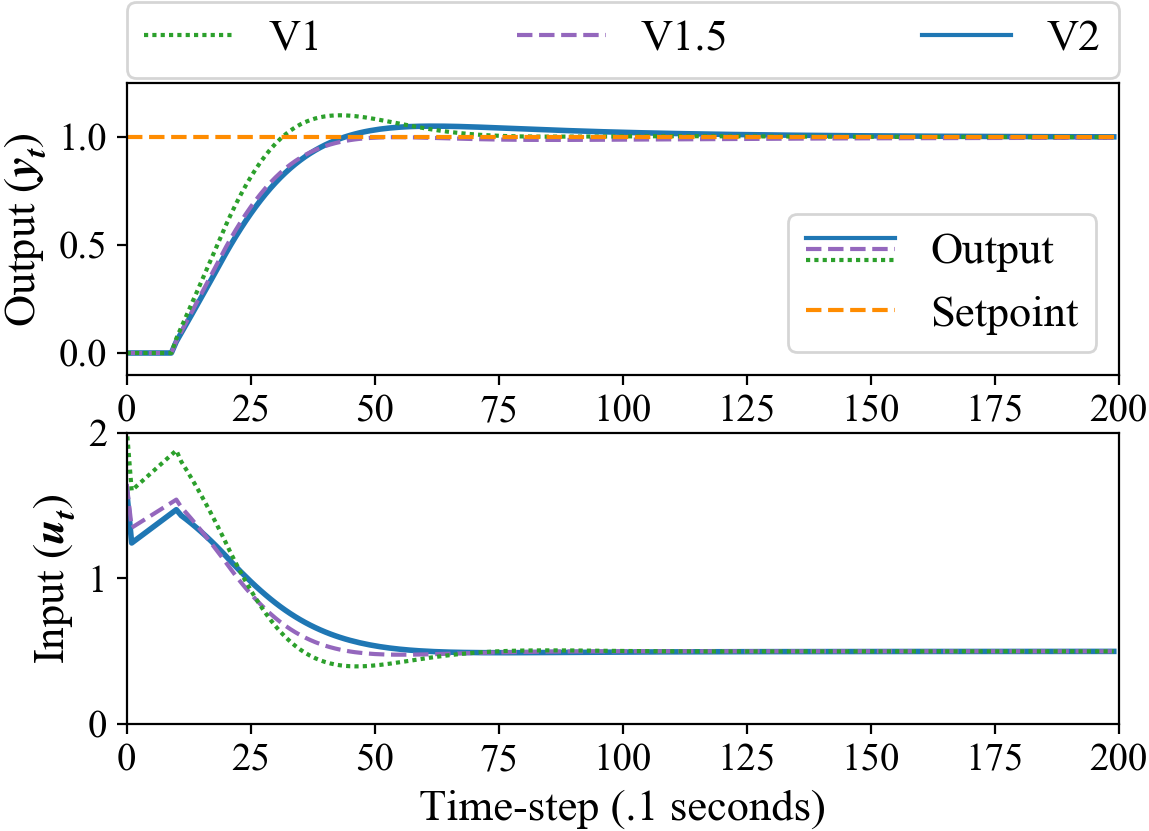}    
\caption{(top) Output signal; (bottom) Input signal. The colors correspond to the respective final PI gains shown in figure \ref{fig:ex1_fig2}.} 
\label{fig:ex1_fig1}
\end{center}
\end{figure}

\begin{figure}
\begin{center}
\includegraphics[width=8.4cm]{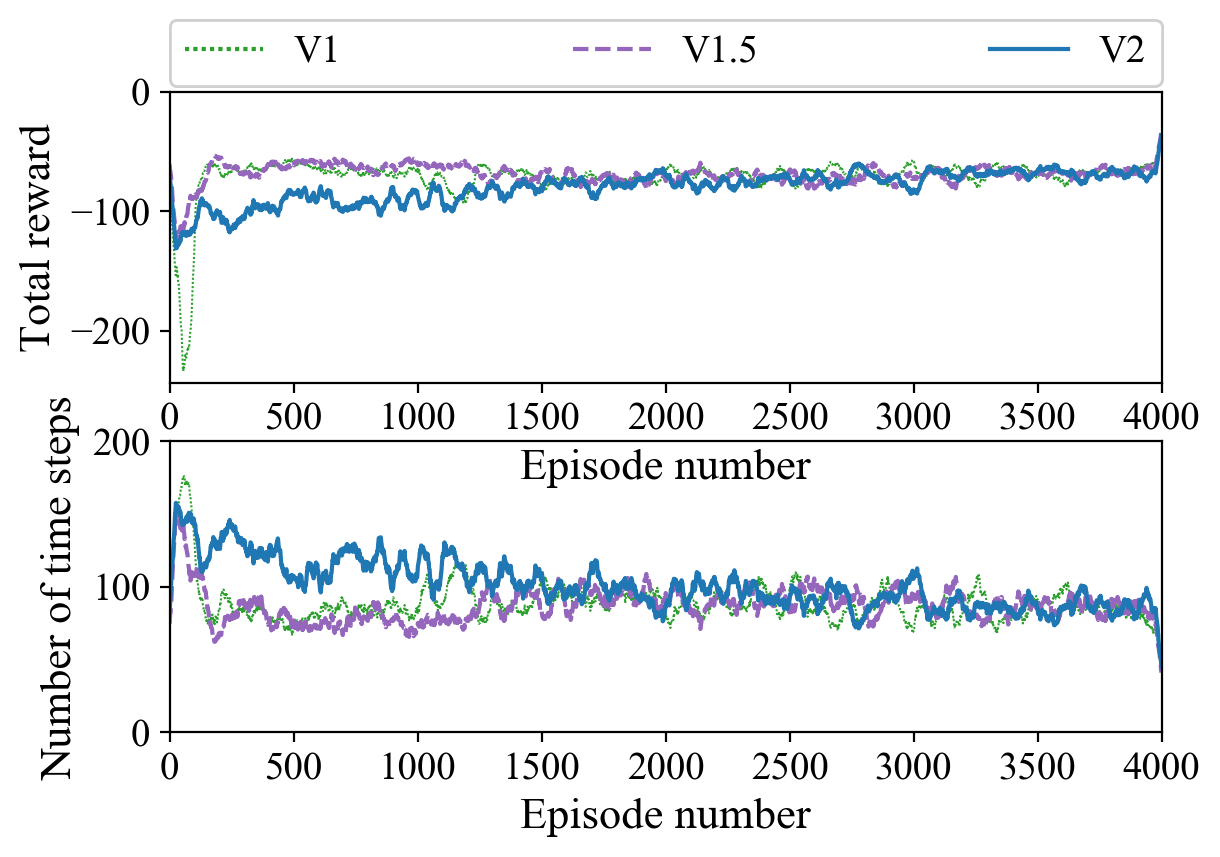}    
\caption{(top) Moving average of total reward per episode; (bottom) Moving average of number of time-steps per episode before the PI controller tracked within 0.1 for 10 consecutive time-steps.} 
\label{fig:ex1_fig4}
\end{center}
\end{figure}

\subsection{Example 1}

Consider the following continuous-time transfer function:
\begin{equation}
    G(s) = \frac{2 e^{-s}}{6s + 1}.
\label{eq:ex1}
\end{equation}
We discretize \eqref{eq:ex1} with time-steps of 0.1 seconds. In this example, we initialize a PI controller with gains $k_p = 0.2, k_i = 0.05$. The following results are representative of other initial PI gains. Note, however, we cannot set $k_p = k_i = 0$ because this forces $\mu( \cdot, K ) \equiv 0$ and the parameters will not change between updates.

In our experiments we implement algorithms \textbf{V1} and \textbf{V2}. We also consider ``\textbf{V1.5}", in which the PID parameters are updated every tenth time-step. Note that the fundamental difference between \textbf{V1} and \textbf{V2} is that the former corresponds to an online implementation of the algorithm, while the latter can be seen as an offline version. \textbf{V1.5} represents a dwell time in the learning algorithm.

Figure \ref{fig:ex1_fig2} only shows the value of $k_p, k_i$ at the end of each episode for each implementation. Nonetheless all three implementations reach approximately the same values; the final closed-loop step responses for each implementation are shown in figure \ref{fig:ex1_fig1}. 

Another way of visualizing the PID parameters is in the $k_p - k_i$ plane. We plot the boundary separating the stable and unstable regions using the parametric curve formulas due to \citet{saeki2007properties}. Figure \ref{fig:ex1_fig5} is a scatter plot with the $k_p, k_i$ value at the end of each episode along with the aforementioned boundary curve. We note that the stability regions refer to the closed-loop with \eqref{eq:ex1} and a fixed $k_p-k_i$ point, rather than the nonlinear system induced by updating $k_p-k_i$ values online.

We see in figure \ref{fig:ex1_fig4} that all three implementations achieve similar levels of performance as measured by the reward function \eqref{eq:reward} ($\lambda = 0.50$). Although, \textbf{V1} and \textbf{V1.5} plateau sooner than \textbf{V2}, the initial dip in reward from \textbf{V1} can be explained by the sequence of unstable $k_p-k_i$ values around the boundary curve in figure \ref{fig:ex1_fig5}.

Finally, we note that \textbf{V1} and \textbf{V1.5} reach their peak performances after approximately 25 minutes (real-time equivalent) of operation. In our experiments, the actor learning rate $\alpha$ had the most drastic effect on the convergence speed. Here, we show the results for a relatively small $\alpha$ (see Appendix \ref{app:details}) to clearly capture the initial upheaval of the parameter updates as well as the long-term settling behavior. In principle, we could omit the latter aspect and simply stop the algorithm, for example, once the reward reaches a certain average threshold.






\begin{figure}
\begin{center}
\includegraphics[width=8.4cm]{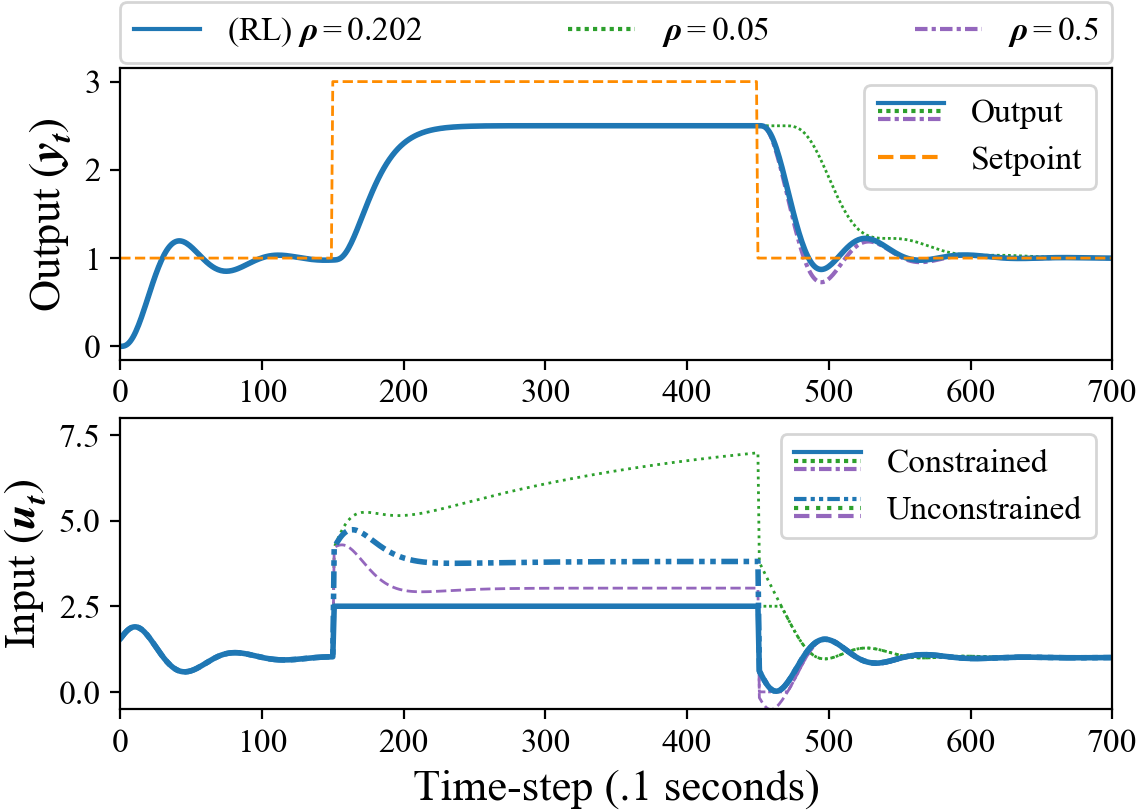}    
\caption{(top) The output response corresponding to various values of $\rho$; (bottom) The colors correspond to the $\rho$ values at the top, while the dashed lines show what the input signal would be without the actuator constraint.} 
\label{fig:ex3_fig1}
\end{center}
\end{figure}

\begin{figure}
\begin{center}
\includegraphics[width=8.4cm]{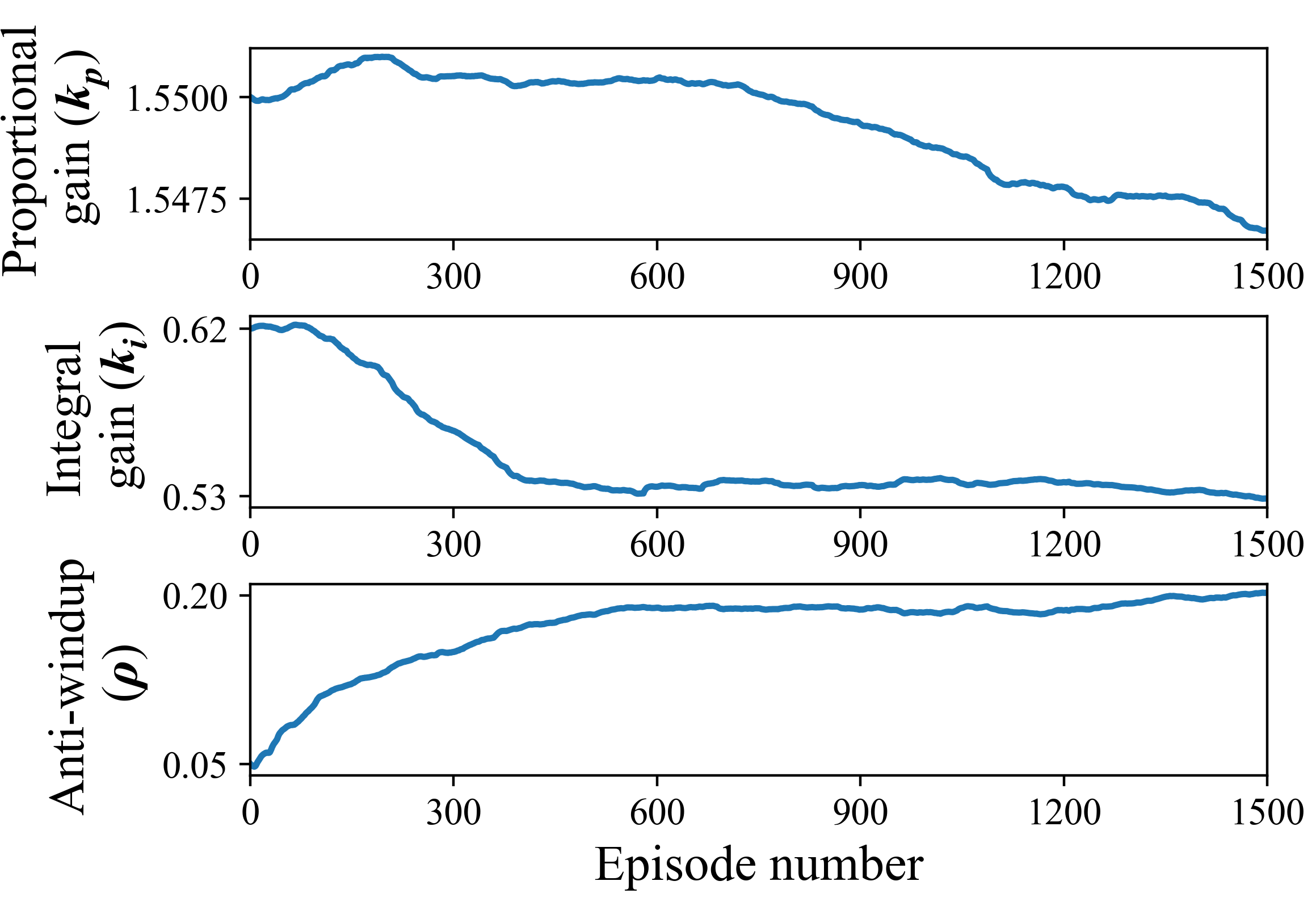}    
\caption{(top) $k_p$ value after each episode; (middle) $k_i$ value; (bottom) $\rho$ value.} 
\label{fig:ex3_fig2}
\end{center}
\end{figure}

\subsection{Example 2}

In this example, we incorporate an anti-windup tuning parameter and employ algorithm \textbf{V2}. Consider the following transfer function:
\begin{equation}
    G(s) = \frac{1}{(s+1)^3}.
\label{eq:ex3}
\end{equation}

In order to tune $\rho$, it is necessary to saturate the input: If the actor always operates within the actuator limits, then
\begin{equation}
\frac{\partial \mu}{\partial \rho} \equiv 0
\label{eq:zeroupdate}
\end{equation}
because $I_u \equiv 0$, meaning $\rho$ will never be updated in \eqref{eq:PolicyGradient_Iteration}. This can be understood from figure \ref{fig:ex3_fig1} at the bottom, as there is a non-zero difference between the dashed and solid lines only after the first step change (this corresponds to the difference shown in figure \ref{fig:Antiwindup}). Further, although \eqref{eq:zeroupdate} also holds for states $s_n$ corresponding to input saturation,
which therefore do not contribute to the update in \eqref{eq:PolicyGradient_Iteration}, we still store them in RM for future policy updates.

In our experiment, the set-point is initialized to $1$, then switches to $1.5$ (plus a small amount of zero-mean Gaussian noise), then switches back to 1. The switches occur at varying time-steps. At the beginning of an episode, with $10\%$ probability, the switches set-point is set to 3 instead of 1.5. Figure \ref{fig:ex3_fig1} shows a slower recovery time for smaller $\rho$ and a more aggressive recovery for larger $\rho$ values. 

We emphasize that the actor can be initialized with hand-picked parameters. To illustrate this, we initialize $k_p$ and $k_i$ using the SIMC tuning rules due to \cite{skogestad2001probably}. Figure \ref{fig:ex3_fig2} shows little change in the $k_p$ parameter, while $k_i$ and $\rho$ adjust significantly, leading to a faster integral reset and smoother tracking than the initial parameters.


\section{Conclusion}

In this work, we relate well-known and simple control strategies to more recent methods in deep reinforcement learning. Our novel synthesis of PID and anti-windup compensation with the actor-critic framework provides a practical and interpretable framework for model-free, DRL-based control design with the goal of being implemented in a production control system. 
Recent works have employed actor-critic methods for process control using ReLU DNNs to express the controller; our work then establishes the simplest, nonlinear, stabilizing architecture for this framework. In particular, any linear control structure with actuator constraints may be used in place of a PID.

\begin{ack}
We would like to thank Profs. Benjamin Recht and Francesco Borrelli of University of California, Berkeley for insightful and stimulating conversations. We would also like to acknowledge the financial support from Natural Sciences and Engineering Research Council of Canada (NSERC) and Honeywell Connected Plant.
\end{ack}

\small
\bibliography{ifacconf}             

\normalsize
\appendix
\section{Implementation Details}\label{app:details}     

In example 1, we use the Adam optimizer to train the actor and critic. To demonstrate simpler optimization methods, we train the actor in example 2 using SGD with momentum (decay constant $0.75$, learning rate decrease $O(1/\sqrt{n})$) and gradient clipping (when magnitude of gradient exceeds 1). Adam, RMSprop, and SGD all led to similar results in all examples. The actor and critic networks were trained using \texttt{TensorFlow} and the processes were simulated in discrete time with the Control Systems Library for \texttt{Python}. The hyperparameters in the DDPG algorithm used across all examples are as follows: Mini-batch size $M = 256$, RM size $10^5$, discount factor $\gamma = 0.99$, initial learning rate for both actor and critic is $0.001$. The critic is modeled by a $64\times64$ ReLU DNN. The saturation function in \eqref{eq:SS3} can be modeled with $\text{ReLU}(x) = \max\{0, x\}$:
$$\text{sat}(u) = \text{ReLU}\big(-\text{ReLU}(u_{\text{max}} - u) + u_{\text{max}} - u_{\text{min}} \big) + u_{\text{min}}.$$

\end{document}